\documentclass[a4paper,reqno]{amsart}

\title{Revisiting multifractal analysis}

\author[Ben Nasr and Peyri\`ere]
{Fathi Ben Nasr$^{*}$ \and Jacques Peyri\`ere$^{\dag\ddag\sharp}$}
\thanks{$^*$ D\'epartement de Math\'ematiques, Facult\'e des Sciences
  de Monastir, Monastir 5000 Tunisie, \texttt{fathi\_bennasr@yahoo.fr}}
\thanks{$^\dag$% Corresponding author\\
Department of Mathematics,
  Tsinghua University, Beijing 100084,
  P.R.~China,\\ \texttt{peyriere@math.tsinghua.edu.cn}}
\thanks{$^\ddag$ Universit\'e Paris-Sud, Math\'ematique b\^at.\ 425,
  CNRS UMR 8628, 91405 Orsay Cedex,
  France. \texttt{jacques.peyriere@math.u-psud.fr}}
\thanks{$\sharp$ Partially suported by Project 111 (P.R.~China)}


\date{}
\usepackage{amsmath,amsfonts,amssymb,mathrsfs}
\newtheorem{theorem}{Theorem}%[section]
\newtheorem{lemma}[theorem]{Lemma}%[section]
\newtheorem{proposition}[theorem]{Proposition}%[section]
\newtheorem{corollary}[theorem]{Corollary}%[section]
\theoremstyle{definition}
\theoremstyle{remark}
\newcommand{\dist}{\mathsf{d}}
\newcommand{\pt}{\tilde{p}}
\newcommand{\entrop}{\mathsf{h}}
\newcommand{\ball}{\mathsf{B}}
\newcommand{\compl}[1]{\vphantom{#1}^\mathsf{C}#1}
\newcommand{\supp}[1]{\mathrm{S}_{#1}} %support
\newcommand{\e}[1]{{\mathrm e}^{#1}}
\DeclareMathOperator*{\esup}{\mathrm{ess\,sup}}
\begin{document}
\begin{abstract}
  New proofs of theorems on the multifractal formalism are given. They
  yield results even at points~$q$ for which Olsen's functions~$b(q)$
  and~$B(q)$ differ. Indeed, we provide an example of measure for
  which functions~$b$ and~$B$ differ and for which the Hausdorff
  dimension of the sets~$E_\alpha$ (the level sets of the local
  H\"older exponent) are given by the Legendre transform of~$b$ and
  their packing dimension by the Legendre transform of~$B$.
\end{abstract}

\maketitle

\section{Introduction}
The multifractal formalism aims at giving expressions of the dimension
of the level sets of the local H\"older exponent of some set
function~$\mu$ in terms of the Legendre transform of some ``free
energy'' function. If such a formula holds, one says that $\mu$
satisfies the multifractal formalism. At first, the formalism used
``boxes'', or in other terms took place in a totally disconnected
metric space. In this context, the closeness to large deviation theory
is patent. To get rid of this boxes and have a formalism meaningful in
Geometric Measure Theory, Olsen~\cite{olsen} introduced a formalism
which nowadays is of common use. At this stage of the theory, whether
it is dealt with or without boxes, the formalism was proven to hold
when there exists an auxiliary measure, a so called \emph{Gibbs
  measure}.  Later on, it was shown that this formalism holds under the
condition that the Olsen's Hausdorff-like multifractal measure be
positive (see~\cite{bn} in the totally disconnected case, \cite{bbh}
in general); So, the situation when $b(q)=B(q)$ (in Olsen's notation)
is fairly well understood.

Here, we elaborate on the previous proofs. There is a vector version
of Olsen's constructs~\cite{peyriere} and in particular of the
functions~$b$ and~$B$; But, this time, they are functions of several
variables. In this work, we show that the restriction to a suitable
affine subspace of these functions allows to estimate the Hausdorff
and Tricot dimensions of some level sets. In particular, this gives
some results even in case when $b\ne B$. Despite the inherent
complexity of notations, not only we provide a simple proof of already
known results, but also we get new estimates.

\section{Notations and definitions}

We deal with a metric space $({\mathbb X},\dist)$ having the
\emph{Besicovitch property:}\\[1ex]
There exists an integer constant~$C_B$ such that one can extract $C_B$
countable families $\bigl\{\{\ball_{j,k}\}_k\bigr\}_{1\le j\le C_B}$
from any collection ${\mathcal B}$ of balls so that
\begin{enumerate}
\item $\displaystyle \bigcup_{j,k} \ball_{j,k}$ contains the centers
  of the elements of~${\mathcal B}$,
\item for any~$j$ and $k\ne k'$, $\ball_{j,k}\cap \ball_{j,k'} =
  \emptyset$.
\end{enumerate}

\subsection*{Notations}

$\ball(x,r)$ stands for the open ball $\ball(x,r) = \{y\in {\mathbb
  X}\ ;\ \dist(x,y)< r\}$. The letter~$\ball$ with or without
subscript will implicitly stand for such a ball.
When dealing with a collection of balls $\{\ball_i\}_{i\in I}$ the
following notation will implicitly be assumed: $\ball_i =
\ball(x_i,r_i)$.

By a $\delta$-cover of~$E\subset {\mathbb X}$, we mean a collection of
\emph{balls} of radii not exceeding~$\delta$ whose union
contains~$E$. A \emph{centered cover} of~$E$ is a cover of~$E$
consisting in balls whose centers belong to~$E$.

By a $\delta$-packing of~$E\subset {\mathbb X}$, we mean a collection
of disjoint balls of radii not exceeding~$\delta$ centered in~$E$.

If~$E$ is a subset of~${\mathbb X}$, $\dim_H E$ stands for its
Hausdorff dimension and~$\dim_P E$ for its packing dimension
(introduced by Tricot~\cite{tricot}).

Let~${\mathscr B}$ stand for the set of balls of~${\mathbb X}$
and~${\mathscr F}$ for the set of maps from~${\mathscr B}$ to
$[0,+\infty)$.

If~$\nu\in {\mathscr F}$, one considers the outer measure~$\nu^{\sharp}$
on~${\mathbb X}$ associated with~$\nu$ in the following way:
\begin{equation*}
\nu^{\sharp}(E) = \inf\left\{\sum \nu(\ball_j)\ ;\ E\subset\bigcup
\ball_j\right\}.
\end{equation*}

The set of $\mu\in {\mathscr F}$ such that~$\mu(\ball)=0$ implies
$\mu(\ball')$ for all~$\ball'\subset \ball$ will be denoted
by~${\mathscr F}^*$. For such a~$\mu$, one defines its
support~$\supp{\mu}$ to be the complement of the set
\begin{equation*}
\bigcup \left\{ \ball\in {\mathscr B}\ ;\ \mu(\ball)= 0 \right\}.
\end{equation*}

\subsection*{Multifractal measures and separator functions}

For~$\mu=(\mu_1,\dots,\mu_m)\in {\mathscr F}^m$, $E\subset {\mathbb
  X}$,~$q=(q_1,\dots,q_m)\in {\mathbb R}^m$, $t\in {\mathbb R}$,
and~$\delta>0$, one sets
\begin{eqnarray*}
\overline{\mathscr P}_{\mu,\delta}^{q,t}(E) &=& \sup \left\{
\sum^{{}\quad *} r_j^t\prod_{k=1}^{m}\mu_k(\ball_j)^{q_k}\ ;\,
\{\ball_j\}\ \delta\text{-packing of }E\right\},\\
\noalign{where~$*$ means that one only sums the terms for which
  $\prod_k\mu_k(\ball_j)\ne 0$,}
\overline{\mathscr P}_{\mu}^{q,t}(E) &=& \lim_{\delta\searrow 0}
\overline{\mathscr P}_{\mu,\delta}^{q,t}(E),\\
{\mathscr P}_{\mu}^{q,t}(E) &=& \inf\left\{ \sum \overline{\mathscr
    P}_{\mu}^{q,t}(E_j)\ ;\ E\subset \bigcup E_j \right\},
\end{eqnarray*}
and% also
\begin{eqnarray*}
\overline{\mathscr H}_{\mu,\delta}^{q,t}(E) &=& \inf \left\{
\sum^{{}\quad *} r_j^t\prod_{k=1}^{m}\mu_k(\ball_j)^{q_k}\ ;\,
\{\ball_j\}\ \text{centered }\delta\text{-cover of }E\right\},\\
\overline{\mathscr H}_{\mu}^{q,t}(E) &=& \lim_{\delta\searrow 0}
\overline{\mathscr H}_{\mu,\delta}^{q,t}(E),\\
%
%\noalign{and}
%
{\mathscr H}_{\mu}^{q,t}(E) &=& \sup\left\{ \overline{\mathscr
  H}_{\mu}^{q,t}(F)\ ;\ F\subset E\right\},
\end{eqnarray*}

When~$m=1$, these measures have been defined by Olsen~\cite{olsen}.
When~$\mu$ is identically~1 these quantities do not depend
on~$q$. They will simply be respectively denoted by
$\overline{{\mathscr P}}_\delta^t(E)$, $\overline{{\mathscr P}}^t(E)$,
${\mathscr P}^t(E)$, $\overline{{\mathscr H}}_\delta^t(E)$,
$\overline{{\mathscr H}}^t(E)$, and ${\mathscr H}^t(E)$. They are the
classical packing pre-measures and measures introduced by
Tricot~\cite{tricot}, and the Hausdorff centered pre-measures and
measures.

Also, as usual, one considers the following functions
\begin{eqnarray*}
\tau_{\mu,E}(q) &=& \inf \{t\in {\mathbb R}\ ;\ \overline{\mathscr
  P}_{\mu}^{q,t}(E)= 0\} = \sup \{t\in {\mathbb
  R}\ ;\ \overline{\mathscr P}_{\mu}^{q,t}(E)= \infty\}\\
B_{\mu,E}(q) &=& \inf \{t\in {\mathbb R}\ ;\ {\mathscr
  P}_{\mu}^{q,t}(E)= 0\} = \sup \{t\in {\mathbb R}\ ;\ {\mathscr
  P}_{\mu}^{q,t}(E)= \infty\},\\
\noalign{and} b_{\mu,E}(q) &=& \inf \{t\in {\mathbb R}\ ;\ {\mathscr
  H}_{\mu}^{q,t}(E)= 0\} = \sup \{t\in {\mathbb R}\ ;\ {\mathscr
  H}_{\mu}^{q,t}(E)= \infty\}.
\end{eqnarray*}

It is well known \cite{olsen,peyriere} that $\tau$ and~$B$ are convex
and that $b\le B\le \tau$. Let $J_\tau$, $J_B$, and $J_b$ stand for
the interior of the sets where respectively~$\tau$,~$B$, and~$b$ are
finite. \medskip

When~$\mu$ is identically~1 we will denote these quantities by
$\dim_B E$, $\dim_P E$, and $\dim_H E$. The first one is the
Minkowski-Bouligand dimension (or box-dimension), the second is the
Tricot (packing) dimension~\cite{tricot}, and the last the Hausdorff
dimension.  \medskip

Here is an alternate definition of~$\tau_{\mu,E}$. Fix $\lambda< 1$
and define

\begin{eqnarray*}
\widetilde{\mathscr P}_{\mu,\delta}^{q,t}(E) &=& \sup \left\{
\sum^{{}\quad *} r_j^t\prod_{k=1}^{m}\mu_k(\ball_j)^{q_k}\,;\,
\{\ball_j\}\ \text{packing of }E \text{ with } \lambda\delta<
r_j\le \delta\right\},\\
\widetilde{\mathscr P}_{\mu}^{q,t}(E) &=& \limsup_{\delta\searrow
  0}\widetilde{\mathscr P}_{\mu,\delta}^{q,t}(E),\\
\noalign{and}
\widetilde{\tau}_{\mu,E}(q) &=& \sup \left\{t\in {\mathbb
  R}\ ;\widetilde{\mathscr P}_{\mu}^{q,t}(E)= +\infty \right\}.
\end{eqnarray*}

\begin{lemma}
One has $\widetilde{\tau}_{\mu,E} = {\tau}_{\mu,E}$.
\end{lemma}

\proof Obviously $\widetilde{\mathscr P}_{\mu}^{q,t}(E)\le
\overline{\mathscr P}_{\mu}^{q,t}(E)$, so $\widetilde{\tau}_{\mu,E}\le
\tau_{\mu,E}$.  To prove the converse inequality, one only has to
consider the case $\tau_{\mu,E}(q)> -\infty$.

Choose $\gamma< \tau_{\mu,E}(q)$ and~$\varepsilon> 0$ such that
$\gamma+\varepsilon < \tau_{\mu,E}(q)$. There exists~$n_0$ such that,
for all~$n>n_0$, there exists a $\lambda^n$-packing $\{\ball_j\}$
of~$E$ such that
\begin{equation*}
\sum r_j^{\gamma+\varepsilon}\prod_{k=1}^{m}\mu_k(\ball_j)^{q_k} > 1.
\end{equation*}
As
\begin{equation*}
\sum r_j^{\gamma+\varepsilon}\prod_{k=1}^{m}\mu_k(\ball_j)^{q_k} =
\sum_{i\ge 0} \sum_{\lambda<r_j\lambda^{-(n+i)}\le 1}
r_j^{\gamma+\varepsilon}\prod_{k=1}^{m}\mu_k(\ball_j)^{q_k},
\end{equation*}
there exists~$i\ge 0$ such that
\begin{equation*}
\sum_{\lambda<r_j\lambda^{-(n+i)}\le 1}
r_j^{\gamma+\varepsilon}\prod_{k=1}^{m}\mu_k(\ball_j)^{q_k} >
\lambda^{i\varepsilon}(1-\lambda^\varepsilon),
\end{equation*}
from which it follows
\begin{equation*}
\sum_{\lambda<r_j\lambda^{-(n+i)}\le 1}
r_j^{\gamma}\prod_{k=1}^{m}\mu_k(\ball_j)^{q_k} >
\lambda^{-(n+i)\varepsilon}\lambda^{i\varepsilon}(1-\lambda^\varepsilon)
= \lambda^{-n}(1-\lambda^\varepsilon),
\end{equation*}
and $\widetilde{\mathscr P}_{\mu}^{q,t}(E)= +\infty$.

\begin{corollary}\label{newdef}
For any $\lambda<1$, one has
\begin{multline*}
  \tau_{\mu,E}(q) = \ \\
  \varlimsup_{\delta\searrow 0}\frac{-1}{\log \delta} \log\sup \left\{
    \sum^{{}\quad *} \prod_{k=1}^{m}\mu_k(\ball_j)^{q_k}\ ;\,
    \{\ball_j\} \text{ packing of }E \text{ with } \lambda\delta<
    r_j\le \delta\right\}.
\end{multline*}
\end{corollary}
\bigskip

\subsection*{Level sets of local H\"older exponents}\ \medskip

Let~$\mu$ be an element of~${\mathscr F}^*$.  For $\alpha,\, \beta\in
\mathbb{R}$, one sets
\begin{eqnarray*}
  \overline{X}_\mu (\alpha) &=& \left\{ x\in \supp{\mu}\ ;\ \limsup_{r\searrow 0}
    \frac{\log \mu\bigl(\ball(x,r)\bigr)}{\log r}\le \alpha \right\},\\
  \underline{X}_\mu (\alpha) &=& \left\{ x\in \supp{\mu}\ ;\ \liminf_{r\searrow 0}
    \frac{\log \mu\bigl(\ball(x,r)\bigr)}{\log r}\ge \alpha \right\},\\
  X_\mu (\alpha,\beta) &=& \underline{X}_\mu (\alpha)\cap
  \overline{X}_\mu (\beta),\\
  \noalign{and}
  X_\mu (\alpha) &=& \underline{X}_\mu (\alpha)\cap \overline{X}_\mu
  (\alpha).
\end{eqnarray*}

\section{Results}

First, one revisits Billingsley and Tricot
lemmas~\cite{billingsley,tricot}.

\begin{lemma}\label{billingsley}
  Let $E$ be a subset of ${\mathbb X}$ and~$\nu$ an element
  of~${\mathscr F}$.
  \begin{enumerate}
  \item[a)] If $B_{\nu,E}(1)\le 0$, then
    \begin{eqnarray}
      \dim_H E &\le& \sup_{x\in E} \liminf_{r\searrow 0} \frac{\log
        \nu\bigl(\ball(x,r)\bigr)}{\log r},\label{bt1}\\
      \dim_P E &\le& \sup_{x\in E} \limsup_{r\searrow 0} \frac{\log
        \nu\bigl(\ball(x,r)\bigr)}{\log r}.\label{bt2}
    \end{eqnarray}
  \item[b)] If $\nu^{\sharp}(E)>0$, then
    \begin{eqnarray}
      \dim_H E &\ge& \esup_{x\in E,\,\nu^\sharp} \liminf_{r\searrow 0}
      \frac{\log \nu\bigl(\ball(x,r)\bigr)}{\log r},\label{bt3}\\
      \dim_P E &\ge& \esup_{x\in E,\,\nu^\sharp} \limsup_{r\searrow 0}
      \frac{\log \nu\bigl(\ball(x,r)\bigr)}{\log r},\label{bt4}
    \end{eqnarray}
where
\begin{equation*}
\esup_{x\in E,\,\nu^\sharp}\chi(x) = \inf \biggl\{t\in{\mathbb R}
;\ \nu^\sharp\Bigl(E\cap\{\chi>t\}\Bigr)=0\biggr\}.
\end{equation*}
  \end{enumerate}
\end{lemma}

\proof Take~$\displaystyle\gamma> \sup_{x\in E} \liminf_{r\searrow 0}
\frac{\log \nu\bigl(\ball(x,r)\bigr)}{\log r}$ and~$\eta>0$. Since
$B_{\nu,E}(1)\le 0$ there exists a partition $\displaystyle E=\bigcup
E_j$ such that $\displaystyle \sum \overline{\mathscr
  P}_{\nu}^{1,\eta/2}(E_j)<1$. It results that $\displaystyle \sum
\overline{\mathscr P}_{\nu}^{1,\eta}(E_j)=0$.

Let~$F$ be a subset of $E_k$ and~$\delta$ a positive number. For
all~$x\in F$, there exists~$r\le \delta$ such
that~$\nu\bigl(\ball(x,r)\bigr)\ge r^\gamma$. By using the Besicovitch
property there exists a centered~$\delta$-cover $\{\ball_{j}\}$ of~$F$,
which can be being decomposed in $C_B$ packings, such that
$\nu(\ball_{j})\ge r_{j}^{\gamma}$. We then have
\begin{equation*}
  \sum r_{j}^{\gamma+\eta}\le \sum r_{j}^{\eta}\nu(\ball_{j})
  \le C_B\overline{{\mathscr P}}_{\nu,\delta}^{1,\eta}(E_k).
\end{equation*}
Therefore we have, $\overline{{\mathscr H}}^{\gamma+\eta}(F)=0$,
${\mathscr H}^{\gamma+\eta}(E_k)=0$, and finally ${\mathscr
  H}^{\gamma+\eta}(E)=0$. Then~\eqref{bt1} easily follows.
\medskip

To prove~\eqref{bt2}, take~$\displaystyle\gamma> \sup_{x\in E}
\limsup_{r\searrow 0} \frac{\log \nu\bigl(\ball(x,r)\bigr)}{\log r}$
and~$\eta>0$. As previously, there exists a partition $\displaystyle
E=\bigcup E_j$ such that $\displaystyle \sum
\overline{\mathscr P}_{\nu}^{1,\eta}(E_j)= 0$.

For all~$x\in E$, there exists~$\delta>0$ such that, for all~$r\le
\delta$, one has~$\nu\bigl(\ball(x,r)\bigr)\ge r^\gamma$. Consider the
set
\begin{equation*}
  E(n) = \left\{ x\in E\ ;\ \forall r\le 1/n,
    \ \nu\bigr(\ball(x,r)\bigr)\ge r^\gamma\right\}.
\end{equation*}

Let $\{\ball_{j}\}$ be a
$\delta$-packing of $E_k\cap E(n)$, with~$\delta\le 1/n$. One has
\begin{equation*}
  \sum r_{j}^{\gamma+\eta}\le
  \sum_{j}r_{j}^\eta\nu(\ball_{j}) \le 
  \overline{{\mathscr P}}_{\nu,\delta}^{1,\eta}(E_k),
\end{equation*}
from which~$\overline{{\mathscr
    P}}^{\gamma+\eta}\bigl(E_k\cap E(n)\bigr)=0$ follows.

So we have ${\mathscr P}^{\gamma+\eta}\bigl(E(n)\bigr)=0$. Since
$\displaystyle E = \bigcup_{n\ge 1} E(n)$, one has $\dim_P E \le
\gamma+\eta$. Hence~\eqref{bt2}.  \medskip

To prove~\eqref{bt3}, take~$\gamma<\esup_{x\in E,\,\nu^\sharp}
\liminf_{r\searrow 0} \frac{\log \nu\bigl(\ball(x,r)\bigr)}{\log r}$
and consider the set $F = \left\{x\in E\ ;\ \liminf_{r\searrow 0}
\frac{\log \nu\bigl(\ball(x,r)\bigr)}{\log r} > \gamma\right\}$. Then, for
all~$x\in F$, there exists~$\delta> 0$ such that, for all~$r\le
\delta$, one has $\nu\bigl(\ball(x,r)\bigr)\le r^\gamma$. Consider the
set
\begin{equation*}
  F(n) = \left\{ x\in F\ ;\ \forall r\le
    1/n,\ \nu\bigr(\ball(x,r)\bigr)\le r^\gamma\right\}.
\end{equation*}
We have $F=\bigcup_{n\ge 1} F(n)$. Since we assume
$\nu^{\sharp}(E)>0$, there exists~$n$ such that $\nu^{\sharp}\bigl(
F(n)\bigr)>0$. Then for any centered $\delta$-cover $\{\ball_j\}$
of~$F(n)$, with $\delta\le 1/n$, one has
\begin{equation*}
  0<\nu^{\sharp}\bigl( F(n)\bigr) \le \sum \nu^{\sharp}(\ball_j) 
  \le \sum \nu(\ball_j)\le   \sum r_j^{\gamma}.
\end{equation*}
Therefore,\quad $\dim_H E\ge \dim_H F(n)\ge \gamma$ (one can compute
the Hausdorff dimension by using centered covers).  \medskip

To prove~\eqref{bt4}, take~$\gamma< \esup_{x\in
  E,\,\nu^\sharp}\limsup_{r\searrow 0} \frac{\log
  \nu\bigl(\ball(x,r)\bigr)}{\log r}$ and consider the set $F =
\left\{x\in E\ ;\ \limsup_{r\searrow 0} \frac{\log
    \nu\bigl(\ball(x,r)\bigr)}{\log r} > \gamma\right\}$. Let~$G$ be a
subset of~$F$. Then, for all~~$x\in G$, for all~$\delta> 0$, there
exists~$r\le \delta$ such that $\nu\bigl(\ball(x,r)\bigr)\le
r^\gamma$. Then for all~$\delta$, by using the Besicovitch property,
there exist a collection $\bigl\{\{\ball_{j,k}\}_j\bigr\}_{1\le k\le
  C_B}$ of $\delta$-packings of~$G$ which together cover~$G$ and
such that $\nu(\ball_{j,k})\le r_{j,k}^{\gamma}$. Then one has
\begin{equation*}
  0<\nu^{\sharp}( G) \le \sum_{j,k} \nu^{\sharp}(\ball_{j,k}) 
  \le \sum_{j,k} \nu(\ball_{j,k})\le
  \sum r_{j,k}^{\gamma}.
\end{equation*}
This implies that there exists~$k$ such that $\displaystyle \sum_j
r_{j,k}^{\gamma}\ge \frac{1}{C_B}\,\nu^{\sharp}(G)$. This implies
$\overline{{\mathscr P}}^\gamma(G) \ge
\frac{1}{C_B}\,\nu^\sharp(G)$. So if $F=\bigcup G_j$, one has
\begin{equation*}
  \sum \overline{\mathscr P}^\gamma(G_j)\ge \frac{1}{C_B}\sum \nu^\sharp(G_j)
  \ge \frac{1}{C_B}\nu^\sharp(F)>0,
\end{equation*}
so ${\mathscr P}^\gamma(F)>0$. Therefore,\quad $\dim_P F\ge \gamma$.
Then~\eqref{bt4} easily follows.

\begin{lemma}\label{main}
  Let~$\mu$ and~$\nu$ be elements of~${\mathscr F}^*$ and~${\mathscr
    F}$ respectively. Set $\varphi(t) = B_{(\mu,\nu),\supp{\mu}}
  (t,1)$ and assume that $\varphi(0)=0$ and $\nu^{\sharp}(\supp{\mu})>
  0$. Then one has
\begin{equation*}
  \nu^{\sharp}\Bigl( \compl{X_\mu
    \bigl(-\varphi'_r(0),-\varphi'_l(0)\bigr)}\Bigr)=0,
\end{equation*}
where $\varphi_l'$ and $\varphi_r'$ are the left and right hand sides
derivatives of~$\varphi$.

The same result holds with $\varphi(t) = \tau_{(\mu,\nu),\supp{\mu}}
(t,1)$.
\end{lemma}

\proof Take $\gamma>-\varphi'_l(0)$, and consider the set
\begin{equation*}
E(\gamma) = \left\{ x\in \supp{\mu}\ ;\ \limsup_{r\searrow 0}
\frac{\log \mu\bigl( \ball(x,r)\bigr)}{\log r} > \gamma\right\}.
\end{equation*}

If~$x\in E(\gamma)$, for all~$\delta>0$, there exists~$r\le \delta$
such that $\mu\bigl( \ball(x,r)\bigr)\le r^\gamma$. Consider a
partition of $E(\gamma)$: \quad $\displaystyle E(\gamma) = \bigcup
E_j$.

For~$\delta>0$, for all~$j$, one can find a $\delta$-cover
$\{\ball_{j,k}\}$ of~$E_j$ such that $\mu(\ball_{j,k})\le
r_{j,k}^{\gamma}$.

We have, for any $t>0$,
\begin{multline*}
  \nu^{\sharp}(E_j) \le \sum_{k} \nu^{\sharp}(\ball_{j,k})
  \le \sum \nu(\ball_{j,k})\\
  = \sum \mu(\ball_{j,k})^{-t} \mu(\ball_{j,k})^t\nu(\ball_{j,k}) \le
  \sum \mu(\ball_{j,k})^{-t} r_{j,k}^{\gamma t}\nu(\ball_{j,k}),
\end{multline*}
which, together with the Besicovitch property, implies
\begin{equation*}
  \nu^{\sharp}\bigl(E(\gamma)\bigr) \le C_B \sum_{j} \overline{\mathscr
    P}_{(\mu,\nu)}^{(-t,1),\gamma t}(E_j)
\end{equation*}
and
\begin{equation*}
  \nu^{\sharp}\bigl(E(\gamma)\bigr) \le C_B 
  {\mathscr P}_{(\mu,\nu)}^{(-t,1),\gamma t}(\supp{\mu}).
\end{equation*}
So, if $\gamma t> \varphi(-t)$, we have
$\nu^{\sharp}\bigr(E(\gamma)\bigr)=0$. But, since $\gamma>
-\varphi'_l(0)$, this happens for small enough positive~$t$.

We conclude that $$\nu^{\sharp}\left( \left\{ x\in \supp \mu\ ;\
    \limsup_{r\searrow 0} \frac{\log \mu\bigl( \ball(x,r)\bigr)}{\log
      r} > -\varphi'_l(0)\right\}\right) = 0.$$

In the same way, one proves that $$\displaystyle \nu^{\sharp}\left(
  \left\{ x\in \supp{\mu}\ ;\ \liminf_{r\searrow 0} \frac{\log
      \mu\bigl( \ball(x,r)\bigr)}{\log r} <
    -\varphi'_r(0)\right\}\right) = 0.$$

\begin{corollary}\label{newthm}
  With the same notations and hypotheses as in the previous lemma, one
  has
\begin{multline*}
  \dim_{H} X_\mu \bigl(-\varphi_r'(0),-\varphi_l'(0)\bigr) \ge\\ \inf
  \left\{ \varliminf_{n\to \infty} \frac{\log
    \nu\bigl(\ball(x,r)\bigr)}{\log r}\ ;\ {x\in
    X\bigl(-\varphi_r'(0),-\varphi_l'(0)\bigr)}\right\}
\end{multline*}
and
\begin{multline*}
\dim_{P} X_\mu \bigl(-\varphi_r'(0),-\varphi_l'(0)\bigr) \ge\\ \inf
\left\{ \varlimsup_{n\to \infty} \frac{\log
  \nu\bigl(\ball(x,r)\bigr)}{\log r}\ ;\ {x\in
  X\bigl(-\varphi_r'(0),-\varphi_l'(0)\bigr)}\right\}.
\end{multline*}
\end{corollary}
\bigskip

%\subsection*{Comparison with the usual formalism}\ \medskip

The previous lemmas contain the nowdays classical results on
multifractal analysis~\cite{olsen,bbh,peyriere}. 

Indeed, let $\mu$ be a element of~${\mathscr F}^*$. Till the end of
this section, we will write~$b$, $\tau$, and $B$ instead of
$b_{\mu,\supp{\mu}}$, $\tau_{\mu,\supp{\mu}}$, and
$B_{\mu,\supp{\mu}}$. For~$q\ge 0$, take
$\nu(\ball)=\mu(\ball)^{q}r^{B(q)}$. Then the corresponding~$\varphi$
of Lemma~\ref{main} is $B_{(\mu,\nu),\supp{\mu}} (t,1) = B(q+t)-B(q)$
and, for~$x\in \overline{X}_\mu(\alpha)$, one has
\begin{equation*}
  \limsup_{r\searrow 0} \frac{\log \nu\bigl(\ball(x,r)\bigr)}{\log r} = 
  q \limsup_{r\searrow 0} \frac{\log \mu\bigl( \ball(x,r)\bigr)}{\log
    r}+B(q)\le q\alpha+B(q).
\end{equation*}
So, due to Lemma~\ref{billingsley}-\eqref{bt2} one gets
\begin{equation*}
\dim_P \overline{X}_{\mu}(\alpha) \le \inf_{q\ge 0} q\alpha+B(q).
\end{equation*}

In the same way, we get

\begin{equation*}
\dim_P \underline{X}_{\mu}(\alpha) \le \inf_{q\le 0} q\alpha+B(q).
\end{equation*}

If moreover we assume that ${\mathscr H}_\mu^{q,B(q)}(\supp{\mu})> 0$,
we have $\nu^{\sharp}(\supp{\mu})>0$, and therefore, due to Lemma~\ref{main}
\begin{equation*}
\nu^{\sharp}\left( \left\{X_\mu\bigl( -B'_r(q),-B'_l(q)\bigr) \right\}\right)
> 0.
\end{equation*}

Therefore, due to Lemma~\ref{billingsley}-\eqref{bt3}, we have
\begin{equation*}
\dim_H \left\{X_\mu\bigl( -B'_r(q),-B'_l(q)\bigr) \right\} \ge
\begin{cases}
-q\,B'_r(q)+B(q)&\text{if } q\ge 0,\\[4pt]-q\,B'_l(q)+B(q)&\text{if }
q\le 0. 
\end{cases}
\end{equation*}

Recall that the Legendre transform of a function~$\chi$ is defined to
be\\ $\displaystyle \chi^*(\alpha) = \inf_{q\in {\mathbb R}}
q\alpha+\chi(q)$.

All this gives a new proof of the following theorem (see \cite{bn} in
the totally disconnected case, \cite{bbh} in general).

\begin{theorem}\label{fbn}
  If $B$ has a derivative at some point~$q\in J_B$ and if ${\mathscr
    H}_{\mu}^{q,B(q)}(\supp{\mu})>0$, then
\begin{equation*}
\dim_H X_\mu \bigl( -B'(q)\bigr) = B^*\bigl(-B'(q)\bigr).
\end{equation*}

The same statement holds with~$\tau$ instead of~$B$.
\end{theorem}

In~\cite{bbh} it is shown that if~$B'(q)$ exists and if $\dim_H X_\mu
\bigl( -B'(q)\bigr) = B^*\bigl( -B'(q)\bigr)$, then $b(q)=B(q)$.
\medskip

We now deal with the case when $b(q)\ne B(q)$. The following notations
will prove convenient: for a real function $\psi$, we set
\begin{equation*}
  \psi_l^\flat(q) = \limsup_{t\searrow 0} \frac{\psi(q-t)-\psi(q)}{-t} 
  \text{\quad and\quad} 
  \psi_r^\flat(q) = \limsup_{t\searrow 0} \frac{\psi(q+t)-\psi(q)}{t}.
\end{equation*}

\begin{lemma}\label{b}
Let~$\mu$ and~$\nu$ be elements of~${\mathscr F}^*$ and~${\mathscr F}$
respectively. Set $\varphi(t) = b_{(\mu,\nu),\supp{\mu}}(t,1)$ and
assume that~$\varphi(0)= 0$ and $\nu^{\sharp}(\supp{\mu})> 0$. Then one has
\begin{eqnarray*}
  &&\nu^{\sharp}\left( \left\{ x\in \supp{\mu}\ ;\ \varliminf_{r\searrow 0}
      \frac{\log \mu\bigl( \ball(x,r)\bigr)}{\log r}> -\varphi_l^\flat(0) 
    \right\}\right) = 0\\
\noalign{and}
&&\nu^{\sharp}\left( \left\{ x\in \supp{\mu}\ ;\ \varlimsup_{r\searrow 0}
    \frac{\log \mu\bigl( \ball(x,r)\bigr)}{\log
      r}< -\varphi_r^\flat(0) \right\}\right) = 0,
\end{eqnarray*}
\end{lemma}

\proof Take $\gamma> -\varphi_l^\flat(0) = \liminf_{t\searrow 0}
\frac{\varphi(-t)}{t}$ and consider the set
\begin{equation*}
E = \left\{ x\in \supp{\mu}\ ;\ \liminf_{r\searrow 0} \frac{\log
  \mu\bigl( \ball(x,r)\bigr)}{\log r}> \gamma \right\}.
\end{equation*}
For all~$x\in E$, there exists~$\delta>0$ such that, for all~$r<
\delta$, one has $\mu\bigl( \ball(x,r)\bigr)< r^\gamma$.

Set $\displaystyle E_{\delta} = \left\{x\in \supp{\mu}\ ;\ \forall
r\le \delta,\ \mu\bigl( \ball(x,r)\bigr)< r^\gamma \right\}$. If
$\{\ball_j\}_j$ is any centered $\delta$-cover of~$E_\delta$, one has,
for any~$t>0$,
\begin{eqnarray*}
\nu^{\sharp}(E_\delta) &\le& \sum \nu^{\sharp}(\ball_j) \le \sum \nu(\ball_j)\\
&\le& \sum\mu(\ball_j)^{-t}\mu(\ball_j)^t\nu(\ball_j)\le
\sum\mu(\ball_j)^{-t}r_j^{\gamma t}\nu(\ball_j)
\end{eqnarray*}
Therefore
\begin{equation*}
\nu^{\sharp}(E_\delta) \le \overline{\mathscr H}_{(\mu,\nu)}^{(-t,1),\gamma
  t}(\supp{\mu})\le {\mathscr H}_{(\mu,\nu)}^{(-t,1),\gamma t}(\supp
\mu).
\end{equation*}
Due to the choice of~$\gamma$ there exists~$t>0$ such that ${\mathscr
  H}_{(\mu,\nu)}^{(-t,1),\gamma t}(\supp{\mu}) = 0$. This proves the
first assertion. The second one is proved in the same way.  \medskip

\begin{proposition}
  Let~$\mu$ be an element in~${\mathcal F}$. Suppose that for
  some~$q\in J_b$, ${\mathscr H}_{\mu}^{q,b(q)}(\supp{\mu})>0$, and
  consider the set
\begin{equation*}
  E = \left\{ x\in\supp{\mu}\ ;\ \varliminf_{r\searrow 0}
    \frac{\log \mu\left( \ball(x,r)\right)}{\log r} \le -b_l^\flat(q)
    \text{ and }\varlimsup_{r\searrow 0}
    \frac{\log \mu\left( \ball(x,r)\right)}{\log r} \ge -b_r^\flat(q)\right\}.
\end{equation*}
Then we have
\begin{equation*}
\dim_P E \ge
\begin{cases}
b(q) - q\,b_r^\flat(q),&\text{if } q\ge 0,\\
b(q) - q\,b_l^\flat(q),&\text{if } q\le 0.
\end{cases}
\end{equation*}
In particular, if $b'(q)$ exists one has
\begin{equation*}
  \dim_P \left\{x\in\supp{\mu} \ ;\  \varliminf_{r\searrow 0}
    \frac{\log \mu\left( \ball(x,r)\right)}{\log r} \le -b'(q) \le
    \varlimsup_{r\searrow 0}
    \frac{\log \mu\left( \ball(x,r)\right)}{\log r}\right\} \ge b(q)-q\,b'(q).
\end{equation*}
\end{proposition}

\proof This results from Lemma~\ref{b} and
Lemma~\ref{billingsley}-\eqref{bt4}

\section{An example}

Now, we can deal with the example given in~\cite{bbh} (Theorem~2.6).We
take for ${\mathbb X}$ the space $\{0,1\}^{{\mathbb N}^*}$ endowed
with the ultrametric which assigns diameter~$2^{-n}$ to cylinders of
order~$n$.

We are given two numbers such that $0< p< \pt\le 1/2$ and a sequence
of integers $1=t_0< t_1< \dots< t_n< \cdots$ such that $\displaystyle
\lim_{n\to \infty} t_n/t_{n+1}= 0$.

We define a probability measure~$\mu$ on $\{0,1\}^{{\mathbb N}^*}$:
the measure assigned to the cylinder
$[\varepsilon_1\varepsilon_2\dots\varepsilon_n]$ is
\begin{equation*}
\mu\bigl( [\varepsilon_1\varepsilon_2\dots\varepsilon_n]\bigr) =
\prod_{j=1}^n \varpi_j,
\end{equation*}
where
\begin{itemize}
\item[-] if $t_{2k-1}\le j< t_{2k}$ for some~$k$, $\varpi_j= p$ if
  $\varepsilon_j=0$, $\varpi_j= 1-p$ otherwise,
\item[-] if $t_{2k}\le j< t_{2k+1}$ for some~$k$, $\varpi_j= \pt$ if
  $\varepsilon_j=0$, $\varpi_j= 1-\pt$ otherwise.
\end{itemize}
\medskip

As a matter of fact, the measure considered in~\cite{bbh} is obtained
by taking the image of $\mu$ under the natural binary coding of
numbers in~$[0,1]$ composed with the Gray code. The purpose of using
the Gray code was to get a doubling measure on~$[0,1]$.
\medskip

For $q\in {\mathbb R}$, define 
\begin{equation*}
  \theta(q) = \log_2 \bigl( p^q+(1-p)^q\bigr) \text{\quad and\quad } 
  \tilde{\theta}(q) = \log_2 \bigl(\pt^q+(1-\pt)^q\bigr).
\end{equation*}
It results from~\cite{bbh} that for $0< q< 1$ we have
\begin{equation*}
b(q)= \theta(q) < \tilde{\theta}(q) = B(q),
\end{equation*}
and, for $q<0$ or~$q>1$,
\begin{equation*}
b(q)= \tilde{\theta(q)} < \theta(q) = B(q).
\end{equation*}

We have the following result.

\begin{proposition}
\begin{enumerate}
\item For $\alpha \in \bigl( -\log_2 (1-\pt),-\log_2 \pt\bigr)$, we
  have
\begin{equation*}
\dim_H X_\mu(\alpha) = \inf_{q\in {\mathbb R}}b(q)+\alpha q.
\end{equation*}
\item For $\alpha \in \bigl( -\log_2 (1-\pt),-\log_2 \pt\bigr)
  \setminus \bigl( [-B'_r(0),-B'_l(0)] \cup
  [-B'_r(1),-B'_l(1)]\bigr)$, we have
\begin{equation*}
\dim_P X_\mu(\alpha) = \inf_{q\in {\mathbb R}}B(q)+\alpha q.
\end{equation*}
\end{enumerate}
\end{proposition}

\proof We consider the measure~$\nu$ constructed as~$\mu$ with
parameters~$r$ and~$\tilde{r}$ instead of~$p$ and~$\pt$. We impose the
condition
\begin{equation}\label{eqderivees}
r\log p + (1-r)\log (1-p) = \tilde{r}\log \pt + (1-\tilde{r})\log
(1-\pt).
\end{equation}
As both~$r$ and~$\tilde{r}$ should belong to the interval~$(0,1)$, we
must have
\begin{equation}\label{var}
  \log \frac{1-p}{1-\pt} < r\log \frac{1-p}{p} < \log \frac{1-p}{\pt}.
\end{equation}

From Corollary~\ref{newdef}, it is easy to compute~$\varphi(x)=
\tau_{(\mu,\nu),\supp{\mu}}$: we have
\begin{equation*}
\varphi(x) = \log_2 \max \left\{\bigl( p^xr+(1-p)^x(1-r)\bigr), \bigl(
\pt^x\tilde{r}+(1-\pt)^x(1-\tilde{r})\bigr) \right\}.
\end{equation*}

Condition~\eqref{eqderivees} implies that~$\varphi'(0)$ exists. We set
\begin{equation}\label{exponent}
\alpha = -\varphi'(0) = -r\log_2 p - (1-r)\log_2 (1-p) = r\log_2
\frac{1-p}{p}- \log_2 (1-p).
\end{equation}
It results from~\eqref{var} that~$\alpha$ can take any value in the
interval $\bigl( -\log_2 (1-\pt),-\log_2 \pt\bigr)$.
\medskip

Besides, the strong law of large numbers shows that we have
\begin{eqnarray*}
\liminf_{n\to \infty} \frac{\log_2 \nu\bigl(\ball(x,2^{-n})\bigr)}{-n}
&=& \min \{\entrop(r) , \entrop(\tilde{r})\}\\ \noalign{and}
\limsup_{n\to \infty} \frac{\log_2 \nu\bigl(\ball(x,2^{-n})\bigr)}{-n}
&=& \max \{\entrop(r) , \entrop(\tilde{r})\}
\end{eqnarray*}
for~$\nu$-almost every~$x$, where we set\quad $\entrop(r) = -r\log_2 r
- (1-r)\log_2 (1-r)$.

Then it results from Lemmas~\ref{main} and~\ref{billingsley}-b that
\begin{eqnarray}
\dim_H X_\mu(\alpha) &\ge& \min \{\entrop(r) ,
\entrop(\tilde{r})\}\label{halfh}\\
\noalign{and}
\dim_P X_\mu(\alpha) &\ge& \max \{\entrop(r) ,
\entrop(\tilde{r})\},\label{halft}
\end{eqnarray}
where~$r$, $\tilde{r}$, and~$\alpha$ are linked by
Relations~\eqref{eqderivees} and~\eqref{exponent}.  \medskip

If~$\alpha$ is defined by~\eqref{exponent}, we have
\begin{equation}\label{alpha}
  \alpha = -\theta'(q)\text{\quad if\quad} q = 
  \frac{\log \frac{1-r}{r}}{\log \frac{1-p}{p}} \text{\quad and\quad}
  \alpha = -\tilde{\theta}'(\tilde{q})\text{\quad if\quad} \tilde{q} = 
  \frac{\log \frac{1-\tilde{r}}{\tilde{r}}}{\log \frac{1-\tilde{p}}{\tilde{p}}}.
\end{equation}

Now, fix~$q$ and~$\tilde{q}$ as above in~\eqref{alpha}. One can check
that, for these values of~$q$ and~$\tilde{q}$, one has
\begin{equation}\label{entropy}
  \theta(q)-q\,\theta'(q) = \entrop(r) 
  \text{\quad and\quad }
  \tilde{\theta}(\tilde{q})-\tilde{q}\,\tilde{\theta}'(\tilde{q}) =
  \entrop(\tilde{r}).
\end{equation}
\medskip

In order to have $\theta(q)=b(q)$, we must have $0< q< 1$, which means
\begin{equation}\label{c1}
  \log_2 \frac{1}{p^p(1-p)^{1-p}} < \alpha < \log_2 \frac{1}{\sqrt{p(1-p)}}. 
\end{equation}
In order to have $\tilde{\theta}(\tilde{q})=b(\tilde{q})$, we must
have $\tilde{q}< 0$ or $\tilde{q}> 1$, which means
\begin{eqnarray}
\alpha &>& \log_2 \frac{1}{\sqrt{\tilde{p}(1-\tilde{p})}}\label{c2}\\
\noalign{or}
\alpha &<& \log_2 \frac{1}{\tilde{p}^{\tilde{p}}(1-\tilde{p})^{1-\tilde{p}}}.
\label{c3}
\end{eqnarray}
One can check that at least one of the conditions~\eqref{c1},
\eqref{c2} and~\eqref{c3} is fulfilled.
\medskip

But for any~$q$ such that $b'(q)$ exists, we have (see~\cite{olsen}
or~\cite{bfp}
\begin{equation}\label{upper}
\dim_H X_\mu\bigl(-b'(q)\bigr) \le b(q)-q\,b'(q).
\end{equation}
The first assertion then results from~\eqref{halfh}, \eqref{entropy},
and~\eqref{upper}.
\bigskip

In order to have $\theta(q)=B(q)$, we must have ${q}< 0$ or
${q}> 1$, which means
\begin{eqnarray*}
  \alpha &>& \log_2 \frac{1}{\sqrt{{p}(1-{p})}} = -B'_l(0)\\
\noalign{or}
\alpha &<& \log_2
\frac{1}{{p}^{{p}}(1-{p})^{1-{p}}} = -B'_r(1).
\end{eqnarray*}
In order to have $\tilde{\theta}(\tilde{q})=B(\tilde{q})$, we must
have $0< \tilde{q}< 1$, which means
\begin{equation*}
  -B'_l(1)=\log_2 \frac{1}{\tilde{p}^{\tilde{p}}(1-\tilde{p})^{1-\tilde{p}}} 
  < \alpha < \log_2
  \frac{1}{\sqrt{\tilde{p}(1-\tilde{p})}} = -B'_r(0).
\end{equation*}
Then Assertion~(2) follows as previously.

\section{The vector case}

As in~\cite{peyriere} instead of $\mu(\ball)^q$ one may consider
expressions of the form $\exp -\langle q,\varkappa(\ball)\rangle$,
where $\varkappa$ takes its values in the dual ${\mathbb E}'$ of a
separable Banach space~${\mathbb E}$ and~$q\in {\mathbb E}$.

Let~$\nu$ be an element of~${\mathscr F}$.  For $E\subset {\mathbb
  X}$,~$q\in {\mathbb E}$, $t\in {\mathbb R}$, and~$\delta>0$, one
sets
\begin{eqnarray*}
  \overline{\mathscr P}_\delta^{q,t}(E) &=& \sup \left\{
    \sum^{{}\quad *} r_j^t\e{-\langle
      q,\varkappa(\ball_j)\rangle}\nu(\ball_j)\ ;\, 
    \{\ball_j\}:\, \delta\text{-packing of }E\right\},\\
\overline{\mathscr P}^{q,t}(E) &=& \lim_{\delta\searrow 0}
\overline{\mathscr P}_\delta^{q,t}(E),\\
{\mathscr P}^{q,t}(E) &=& \inf\left\{ \sum \overline{\mathscr
    P}^{q,t}(E_j)\ ;\ E\subset \bigcup E_j \right\},
\end{eqnarray*}
and% also
\begin{eqnarray*}
  \overline{\mathscr H}_\delta^{q,t}(E) &=& \inf \left\{
    \sum^{{}\quad *} r_j^t\e{-\langle
      q,\varkappa(\ball_j)\rangle}\nu(\ball_j)\ ;\
    \{\ball_j\}\ \text{centered }\delta\text{-cover of }E\right\},\\
\overline{\mathscr H}^{q,t}(E) &=& \lim_{\delta\searrow 0}
\overline{\mathscr H}_\delta^{q,t}(E),\\
{\mathscr H}^{q,t}(E) &=& \sup\left\{ \overline{\mathscr
    H}^{q,t}(F)\ ;\ F\subset E\right\},
\end{eqnarray*}

For a function~$\chi$ from ${\mathbb E}$ to~${\mathbb R}$, and
for~$v\in {\mathbb E}$ of norm~1, one defines
\begin{eqnarray*}
\partial_v\chi(0) &=& \lim_{t\searrow 0}\frac{\chi(tv)-\chi(0)}{t}\\
\noalign{and}
\partial_v^*\chi(0) &=& \limsup_{t\searrow 0}-\frac{\chi(tv)-\chi(0)}{t}.
\end{eqnarray*}

With these notations we have the following analogues of
Lemmas~\ref{main} and~\ref{b}.

\begin{lemma}
Let~$\varphi(q)$ be one of the following functions:
$$
  \inf \left\{ t\ ;\ \overline{\mathscr P}^{q,t}({\mathbb X})=0\right\}\\
  \text{\quad or\quad}
  \inf \left\{ t\ ;\ {\mathscr P}^{q,t}({\mathbb X})=0\right\}.
$$
Assume that~$\varphi(0)=0$ and that $\partial_v\varphi(0)$ at~$0$ is a
lower semi-continuous function of~$v$. Then one has
\begin{equation*}
  \nu^{\sharp} \left\{ x\ ;\ \liminf_{r\searrow 0}\frac{\langle
      v,\varkappa\bigl( \ball(x,r)\bigr)}{-\ln r}< -\partial_v\varphi(0)
    \text{\ for some\ } v\in {\mathbb E}  \right\} = 0.
\end{equation*}
\end{lemma}

\begin{lemma}
  Set $\varphi(q) = \inf \left\{ t\ ;\ {\mathscr H}^{q,t}({\mathbb
      X})=0\right\}$ and assume that $\varphi(0)=0$ and that
  $\partial_v^*\chi(0)$ is a lower semi-continuous function
  of~$v$. Then one has
\begin{equation*}
  \nu^{\sharp} \left\{ x\ ;\ \limsup_{r\searrow 0}\frac{\langle
      v,\varkappa\bigl( \ball(x,r)\bigr)}{-\ln r}< -\partial_v^*\varphi(0)
    \text{\ for some\ } v\in {\mathbb E}  \right\} = 0.
\end{equation*}
\end{lemma}

The proofs follow the same lines as above and as the proofs
in~\cite{peyriere}. As a corollary we get the following result (with
the notations of~\cite{peyriere}).

\begin{theorem}
  Let $B(q) = \inf \{t\in {\mathbb R}\ ;\ {\mathscr
    H}_{\varkappa}^{q,t} ({\mathbb X})= 0\}$. Assume that, at some
  point~$q$, the function~$B$ is differentiable with
  derivative~$B'(q)$ and that ${\mathscr
    H}_\varkappa^{q,B(q)}({\mathbb X})> 0$. Then one has
  \begin{equation*}
    \dim_{H} \left\{x \ ;\ \forall v\in {\mathbb E}, 
      \lim_{r\searrow 0} \frac{\bigl\langle v,
        \varkappa\bigl(\ball(x,r)\bigr)\bigr\rangle}
      {\log r} = -B'(q)v\right\} = B(q)-B'(q)q.
  \end{equation*}
\end{theorem}

\end{document}